\documentclass{article}
\begin{document}
\centerline{\textbf{\Large Note on a Nonlinear Differential Equation}}
\[
\]
\centerline{\bf Nikos Bagis}
\centerline{\bf nikosbagis@hotmail.gr}
\[
\]
\centerline{\textbf{Abstract}}
We give evaluations in closed form of certain non linear differential equations 
\[
\]

\section{The equation $y'''(x)+y'(x)=P(y(x))$}

We consider the folowing non linear differential equation (NLDE) 
\begin{equation}
y'''(x)+y'(x)=\frac{1}{2}P(y(x))
\end{equation}
where $P(x)$ is polynomial or a function of the form
\begin{equation} 
P(x)=\left(\sum^{N}_{n=0}a_nx^n\right)\sqrt{\sum^{M}_{n=0}b_nx^n}
\end{equation}
The non linear equation (1) is equivalent to
\begin{equation}
h''(x)=\frac{P(x)}{\sqrt{h(x)-x^2}}.
\end{equation}
and if the above equation has algebraic solution 
$$
h(x)=x^2+U^2(x),
$$ 
then 
$$
y_i(x)=\frac{1}{2}\int\frac{h''(x)}{P(x)}dx
$$
where $y_i$ denotes the inverse of $y(x)$ i.e $y_i(x)=y^{(-1)}(x)$.\\
\\
\textbf{Proof.}\\We define by $F(x)$ the function such that 
\begin{equation}
y''(x)=F(y(x))-y(x)
\end{equation}
then (4) has solution 
\begin{equation}
\int^{x}_{C_3}\frac{1}{\sqrt{C_1+2\int^{t}_{C_2}(F(v)-v)dv}}dt=y_i(x)
\end{equation}
Also relating (1) with (5) we get  
\begin{equation}
\int^{x}_{C}\frac{F'(t)}{P(t)}dt=yi(x)
\end{equation}
Again relating (4) with (5) we have 
\begin{equation}
\frac{F'(x)}{P(x)}=\frac{1}{\sqrt{C_1+2\int^{x}_{C_2}(F(v)-v)dv}}
\end{equation}
Then if
\begin{equation}
h(x)-x^2=U^2(x)=C_1+2\int^{x}_{C_2}(F(t)-t)dt
\end{equation}
we get
$$
h''(x)=\frac{P(x)}{\sqrt{h(x)-x^2}}
$$
\\
The idea behind the above formulation is that last equation admits polynomial solutions if $P(x)$ is of type (2).\\
Let
\begin{equation}
P(x)=
\left(\sum^{N}_{n=0}a_nx^n\right)\sqrt{\sum^{M}_{n=0}b_nx^n},
\end{equation}   
then $b_n$ are depent from $a_n$ with very simple way.\\ 
Set where $h(x)$ in (3) the polynomial $x^2+U(x)^2$, with 
\begin{equation}
U(x)=\sqrt{\sum^{N}_{n=0}b_nx^n}
\end{equation}   
and evaluate the $a_n$ such 
\begin{equation}
h''(x)-\sum^{N}_{n=0}a_nx^n=0
\end{equation}
or equivalently
\begin{equation}
2+2U{'}^2(x)+2U(x)U''(x)=\sum^{N}_{n=0}a_nx^n
\end{equation}
since
\begin{equation}
h(x)=x^2+U^2(x) .
\end{equation}
This can be done as in the following\\
\\
\textbf{Example 1.} Set\\
$P(x)=(g+fx+ex^2)(g_1+f_1x+e_1x^2)$ and $U(x)=g+fx+ex^2$, then if $h(x)=x^2+U^2(x)$, $e_1=12e^2$, $f_1=12ef$, $g_1=2+2f^2+4eg$, we get 
$h(x)=x^2+(g+fx+ex^2)^2$.\\Hence if 
we have to solve the equation
$$
y'''(x)+y'(x)=g+2eg^2+(e+8e^2g)y^2(x)+6e^3y^4(x)
$$
then $P(x)=g+f^2g+2eg^2+(f+f^3+8efg)x+(e+7ef^2+8e^2g)x^2+12e^2fx^3+6e^3x^4$, $h(x)=x^2+(g+fx+ex^2)^2$ and 
$$
y_i(x)=\int\frac{h''(x)}{P(x)}dx=\frac{2 \arctan \left(\frac{2 e x+f}{\sqrt{4 e g-f^2}}\right)}{\sqrt{4 e g-f^2}}
$$
Finaly inverting we get
$$
y(x)=\frac{\sqrt{4eg-f^2}\tan\left(\frac{x}{2}\sqrt{4e g-f^2}\right)-f}{2e}
$$
\\
\textbf{Example 2.} 
The equation
$$
y'''(x)+y'(x)=dy^3(x)+15d^3 y^7(x)
$$
have $P(x)=2dx^3+15d^3x^7$ and $h(x)=x^2+d^2x^6$
and hence 
$$
y(x)=\frac{i}{\sqrt{2d\cdot x}}.
$$ 
Also another equation is with
$$
P(x)=15d^3x^7+21d^2g x^4+dx^3+6dg^2 x+g
$$
then $h(x)=x^2+(g+dx^3)^2$ and the solution is such that
$$
y_i(x)=
$$
$$
=\frac{1}{3 \sqrt[3]{d} g^{2/3}}[45 d^{4/3} g^{2/3} x^2-\left(9 d^{2/3} g^{4/3}+1\right) \log \left(d^{2/3} x^2-\sqrt[3]{d} \sqrt[3]{g} x+g^{2/3}\right)+
$$
$$
+2 \left(9 d^{2/3} g^{4/3}+1\right) \log \left(\sqrt[3]{d} x+\sqrt[3]{g}\right)+2 \sqrt{3} \left(9 d^{2/3} g^{4/3}-1\right) \arctan\left(\frac{1-\frac{2 \sqrt[3]{d} x}{\sqrt[3]{g}}}{\sqrt{3}}\right)]
$$
\\
\textbf{Example 3.}
Assume the NLDE with
$$
P(x)=\frac{1}{2 \sqrt{210}}\left(a_1 x^6+b_1 x^5+c_1 x^4+d_1 x^3+e_1 x^2+f_1 x+g_1\right)\times
$$
$$ 
\times[15 a_1 x^8+20 b_1 x^7+28 c_1 x^6+42 d_1 x^5+70 e_1 x^4+140 f_1 x^3+420 g_1 x^2-
$$
$$
-840 x^2+840 C_2 x+1680 C_1]^{1/2}
$$
then
$$
C+x=2 \sqrt{210}\int^{y(x)}_{0}[15 a_1 t^8+20 b_1 t^7+28 c_1 t^6+42 d_1 t^5+
$$
$$
+70 e_1 t^4+140 f_1 t^3+(420 g_1-840) t^2+840 C_2 t+1680 C_1]^{-1/2}dt
$$
\\
An application of the above example is taking
$$
P(x)=\frac{(f_1 x+g_1) \sqrt{1680 C_1+840 C_2 x+140 f_1 x^3+420 g_1 x^2-840 x^2}}{2 \sqrt{210}}
$$
then
$$
h(x)=x^2
+\frac{1}{840} \left(1680 C_1+840 C_2 x+(420 g_1-840) x^2+140f_1 x^3\right)
$$
and $y(x)$ is given from
$$
x+C=\int^{y(x)}_{0}\frac{dt}{\sqrt{\frac{f_1 }{6}t^3+\frac{(g_1-2)}{2}t^2+C_2t+2C_1}}
$$
The above integral can evaluated using the incomplete elliptic integral of the first kind $F[x,m]$, i.e 
\begin{equation}
F[x,m]=\int^{x}_{0}\frac{dt}{\sqrt{1-m\sin^2(t)}}.
\end{equation}
The reader can see [1].

\section{References}

[1]: J.V. Armitage W.F. Eberlein. 'Elliptic Functions'. Cambridge University Press. (2006)

\end{document}